\documentclass[10pt]{article}
\usepackage[margin=1in]{geometry}
\usepackage{subfiles}
\usepackage[]{mdframed}
\usepackage{float}

\usepackage{amsmath}
\usepackage{amsfonts}
\usepackage{amssymb}
\usepackage{amsthm}
\usepackage{mathtools}

\usepackage{comment}
\usepackage{caption}
\usepackage{subcaption}

\theoremstyle{definition}

\newtheorem*{acknowledgement}{Acknowledgement}

\newcommand{\figref}[1]{\figurename~\ref{#1}}

\usepackage[
backend=biber,
style=alphabetic,
sorting=ynt
]{biblatex}
\begin{document}

\title{A Singular Integral for a Simplified Clairaut Equation}
\author{
  Anand Ganesh \thanks{National Institute of Advanced Studies, email: anandg@nias.res.in} \and
  Anand Rajagopalan \thanks{Cruise, email: anandbr@gmail.com}
}

\maketitle

\begin{abstract}

This expository article on the Lagrange singular integral contains two novelties. The first novelty involves a connection between the Lagrange singular integral for a simplified Clairaut equation, and Euler's homogeneous function theorem. The paper presents a formal derivation of Euler's solution from Lagrange's complete integral, though with some caveats, and also constructs more general surfaces from the complete integral which go beyond Euler's solutions. The first rather complicated construction is based directly on Goursat's definition of a general integral, while the subsequent simpler constructions are based on a suitably expanded notion of the general integral. This generalized general integral is our second novelty. It bridges some of the gap between the the general integral, and the complete integral, partially addressing Evans' remarks (Partial Differential Equations, AMS Graduate Studies in Mathematics, 1998) on the limitations of the general integral. Finally we discuss some subtleties around complete integrals as noted by Chojnacki (Proceedings of the AMS, 1995) and some around general integrals as noted by Evans, and how they apply to our examples. We aim to present these classical PDE concepts to readers with a basic knowledge of multivariable calculus.
\end{abstract}

\section{Introduction}

We look at general solutions to the following PDE (Partial Differential Equation):
\begin{align} \label{eq:restricted_clairaut}
  x\frac{\partial z}{\partial x} + y\frac{\partial z}{\partial y} = z.
\end{align}

One family of solutions to this equation is provided by Euler's homogeneous function theorem. The theorem starts with a function $f(x, y)$ such that $f(sx, sy) = s^nf(x, y)$, in other words, a homogenous function of degree $n$. The theorem states that for such a function,

\begin{align} \label{eq:euler_equation}
  x\frac{\partial f}{\partial x} + y\frac{\partial f}{\partial y} = nf(x, y).
\end{align}

More specifically, if we set $n=1$ we have the following:
\begin{align} \label{eq:euler_equation_1}
  x\frac{\partial f}{\partial x} + y\frac{\partial f}{\partial y} = f(x, y).
\end{align}

This looks remarkably similar to the simplified Clairaut equation \eqref{eq:restricted_clairaut} where $z$ has been replaced by $f(x, y)$. It is easy to verify that functions $z = f(x, y)$ where $f$ is of degree $n=1$ are valid solutions of equation \eqref{eq:restricted_clairaut}. We can also state a partial converse that the most general solution to equation \eqref{eq:euler_equation} is a homogenous function $f$ of degree $n$, or in the case of \eqref{eq:euler_equation_1}, a homogeneous function $f$ of degree $n=1$.

We note an important distinction here, that the two equations \eqref{eq:restricted_clairaut} and \eqref{eq:euler_equation_1} look similiar, but they are not the same, and are not equivalent. In particular, we may have relations $f(x, y, z) = 0$ that are not representable as a {\em function} $z = f(x, y)$ and still obey the restricted Clairaut equation \eqref{eq:restricted_clairaut}, for instance the surface $z^2 = x^2 + y^2$, a double cone. In this paper we refer to such a relation suggestively as a homogeneous surface since the associated polynomial, in this case $z^2 - x^2 - y^2$, is homogeneous. The word {\em surface} is used in contrast to the word {\em function}, and may include singular surfaces like the cone depending on the context.

Another source of confusion is Euler's principle \cite{Engelsman1980} that the integration of a PDE of order $n$ is complete if the integral contains $n$ arbitrary functions. One might argue that the homogenous solution $z=f(x, y)$ contains such an arbitrary function $f$, and thus assume that it is the most general solution possible. This is not true based on various examples we show including the double cone above, $z^2 = x^2 + y^2$. It is possible that homogeneous functions cannot be regarded as arbitrary because of their homogeneity, but it is still unclear how to apply the principle suitably in our context, and whether there is a version that applies to relations $f(x, y, z) = 0$ and not just to functions $z = f(x, y)$.

In order to introduce Lagrange's solution, we start with the most obvious solution to \eqref{eq:restricted_clairaut}, namely $z = ax + by$ with arbitrary constants $a$, $b$. A little more thought may suggest some more solutions like $\sqrt{xy}$ or even $x^{\alpha}y^{1-\alpha}$ all of which are within the ambit of Euler's homogeneous solution. The surprising fact is that the simplest of these solutions, namely $z=ax + by$, can be used to generate all the other solutions including Euler's in an algorithmic manner, and this is the gist of Lagrange's solution.

The current paper makes two claims, one regarding Lagrange's complete integral, and the other regarding Goursat's general integral method which is a specific procedure related to the complete integral. The first claim is that Lagrange's complete integral can generate Euler's homogenous functions, though with some caveats, and also more general surfaces that go beyond these homogenous functions. The second claim regarding the general integral method is that it can be extended suitably to cover a broader range of solutions as described by the complete integral. This addresses some remarks by Evans \cite[p.96]{Evans1998} on the limitations of the general integral. Further, \cite{Chojnacki1995} notes that certain questions about the complete integral have not been well studied, and the current paper fills some gaps in this area.

While the Lagrange integral can formally generate Euler's homogenous functions, there are examples where the envelope procedure breaks down in subtle ways, restricting the domain of applicability. Thus one cannot claim that the Lagrange integral is strictly more general than Euler's homogenous function despite formal generality. These issues reflect Chojnacki's point \cite{Chojnacki1995} that it is difficult to make general statements. Going beyond Euler's solutions, Lagrange's complete integral can be used to generate surfaces like the following one which are not in Euler's functional form $z=f(x, y)$:
\begin{gather*}
  x^2 + y^2 + z^2 - 2xz - 2yz = 0
\end{gather*}

The equation at hand \eqref{eq:restricted_clairaut} is a special case of the Clairaut equation (\cite[p.93]{Evans1998}) which is as follows:
\begin{align} \label{eq:general_clairaut}
  x\frac{\partial z}{\partial x} + y\frac{\partial z}{\partial y} + k\left( \frac{\partial z}{\partial x}, \frac{\partial z}{\partial y} \right) = z.
\end{align}

It is one of the simpler nonlinear PDEs one might see early on in a PDE book. The restricted Clairaut equation provides a good pedagogical platform to discuss classical PDE ideas like the complete integral, the general integral and the singular integral. While Clairaut was a reputed physicist, the motivation for this particular equation appears to have been purely mathematical.

The rest of the paper is organized as follows. The first four sections are largely expository, defining terms related to the singular integral, and describing procedures for constructing them. Section \ref{sec:euler_lagrange_connection} contains both the new claims, the first connecting Lagrange's complete integral to Euler's homogenous function, and the second describing  extensions to the general integral. The final section relates these to the comments of Chojnacki and Evans.

\section{Preliminaries}

To understand Lagrange's solution we need a couple of notions which we now describe. First is the notion of a {\em complete integral}, which an older text like Goursat defines as a family of integrals $f(x, y, z, a, b) = 0$ that involve two arbitrary constants (or parameters) $a$, $b$. Newer texts like \cite[p.92]{Evans1998} have an equivalent definition where the rank of a certain matrix captures the number of independent parameters and equals two in this case. Next we have the notion of a {\em singular integral} which \cite[p.95]{Evans1998} defines as the envelope of this family of integrals (also \cite[p.237]{Goursat1917}). A third notion is that of a {\em general integral} which is a type of complete integral where $b$ is assumed to be a function of $a$, say $b = \phi(a)$ \cite[p.238]{Goursat1917}.

The complete integral for equation \eqref{eq:general_clairaut} is known to be $z = ax + by + k(a, b)$ where $a$ and $b$ are arbitrary constants. For the restricted equation with $k=0$, the complete integral is more specifically:
\begin{align} \label{eq:complete_integral}
  z = ax + by.
\end{align}

According to Lagrange (and \cite[p.236]{Goursat1917}), once we have the complete integral, that is, a family of integrals which depends upon two arbitrary parameters $a$ and $b$, we can derive all the other integrals from them by differentiations and eliminations. But it is not clear how to operationalize this idea, and come up with the solutions indicated earlier. It is also unclear how Euler's homogenous solution relates to this complete integral. The current expository paper is an attempt to describe Goursat's method of differentiations and eliminations, and to connect these integrals to Euler's homogenous solution.

The current problem illustrates two different ways of representing a general solution. According to Euler \cite{Engelsman1980}, an integral of order $n=1$ is complete if it contains $n=1$ arbitrary functions. In his view, the arbitrary constants in the solution to an ordinary differential equation is taken over by an arbitrary function in the case of partial differential equations. Lagrange on the other hand uses the modern notion of a complete integral requiring two arbitrary constants, say $a$, $b$ and a relation $f(x, y, z, a, b) = 0$.

Euler's homogeneous function theorem seems innocuous, and does not claim to present a general solution to a PDE. On the other hand, it is unclear how to apply Euler's principle of arbitrary functions in our context. A careless application of the principle would suggest that homogeneous solutions are the most general ones for the restricted Clairaut equation. But this conclusion is incorrect, as we have seen. As we note in Section \ref{sec:evans_comments}, Evans' use of Euler's principle is also problematic. Further study is required to understand what Euler had in mind, and whether in fact, the principle is valid.

In any case, Lagrange's complete integral seems like an elegant improvement over Euler's principle of arbitrary functions, and does not suffer from such ambiguities. The Lagrange viewpoint seems well matched with modern ideas in differential geometry, and it is not surprising that it continues to find a place in contemporary PDE textbooks. Though it is less ambiguous, it is not well studied whether any given solution to a first order PDE can be derived from Lagrange's complete integral \cite{Chojnacki1995}.

The current paper claims that Lagrange's formulation is formally more general than Euler's homogenous function, though with some important caveats when all the partial derivatives vanish. Euler's solutions are of the form $z = f(x, y)$ where $f$ is a homogenous function of degree $1$. The formulation by Lagrange says that all solutions can be found from the complete integral \eqref{eq:complete_integral} by differentiations and eliminations. In particular, a solution to the equation is not required to be of the form $z = f(x, y)$ where $z$ is a {\em function} of $(x, y)$. In fact, we will find that Lagrange's solution, unlike Euler's, includes a broader class relations $f(x, y, z) = 0$ which may have a one to many relation between $(x, y)$ and $z$.

\section{Euler's Solution} \label{sec:euler_rajaram}
In this section, we will look at a simple and systematic approach to deriving Euler's solution. We start with a change of variables $X = \ln x$, $Y = \ln y$ and $Z = \ln z$. This can be motivated by the suggestive form $\frac{\partial z}{z} \frac{x}{\partial x} + \frac{\partial z}{a} \frac{y}{\partial y} = 1$. With this change of variables, we have the following calculations:
\begin{gather*}
  \frac{\partial Z}{\partial X} = \frac{\partial Z}{\partial z} \frac{\partial z}{\partial x} \frac{\partial x}{\partial X} = \frac{1}{z} \frac{\partial z}{\partial x} x, \\
  \frac{\partial Z}{\partial Y} = \frac{\partial Z}{\partial z} \frac{\partial z}{\partial y} \frac{\partial y}{\partial Y} = \frac{1}{z} \frac{\partial z}{\partial y} y.
\end{gather*}

This gives us the following linearized form of the restricted Clairaut equation:
\begin{gather*}
  \frac{\partial Z}{\partial X} + \frac{\partial Z}{\partial Y} = \frac{x}{z} \frac{\partial z}{\partial x} + \frac{y}{z} \frac{\partial z}{\partial y} = 1.
\end{gather*}

The above logarithmic substitutions restrict the domain of validity to be $x > 0, y > 0, z > 0$. A further change of variables is as follows:
\begin{gather*}
  X = U + V, \\
  Y = U - V.
\end{gather*}

This gives us
\begin{align*}
  \frac{\partial Z}{\partial U} &= \frac{\partial Z}{\partial X}\frac{\partial X}{\partial U} + \frac{\partial Z}{\partial Y}\frac{\partial Y}{\partial U} \\
  &= \frac{\partial Z}{\partial X} + \frac{\partial Z}{\partial Y} \\
  &= 1. \\
\end{align*}

On the other hand,
\begin{align*}
  \frac{\partial Z}{\partial V} &= \frac{\partial Z}{\partial X}\frac{\partial X}{\partial V} + \frac{\partial Z}{\partial Y}\frac{\partial Y}{\partial V} \\
  &= \frac{\partial Z}{\partial X} - \frac{\partial Z}{\partial Y}.
\end{align*}

This is a generic identity and does not seem to impose a particular relationship between $Z$ and $V$. Put differently, we can let $Z$ be an arbitrary function of $V$, and the above identity will continue to hold. Further, $Z=h(V)$ satisfies the following equation for $V = \frac{X-Y}{2}$ and $h$, an arbitrary function:
\begin{gather*}
  \frac{\partial h(V)}{\partial X} + \frac{\partial h(V)}{\partial Y} = \frac{\partial h(V)}{\partial V}\frac{\partial V}{\partial X} + \frac{\partial h(V)}{\partial V}\frac{\partial V}{\partial Y} = h'(V) \frac{1}{2} + h'(V) \frac{-1}{2} = 0.
\end{gather*}

Now $\frac{\partial Z}{\partial U} = 1$ has solutions of the form $Z = U + c$ where $c$ is independent of $U$. Combining this with an arbitrary function $h(V)$, we have the following general solution:
\begin{gather} \label{eq:rajaram_general_uv}
  Z = U + h(V).
\end{gather}

The general solution \eqref{eq:rajaram_general_uv} can also be written as:
\begin{gather*} \label{eq:rajaram_general_xy}
  \ln z = \frac{\ln x + \ln y}{2} + h(\frac{\ln x - \ln y}{2}), \\
  z = \sqrt{xy} \cdot H(\frac{x}{y}).
\end{gather*}

Given that $H$ is arbitrary, other equivalent forms of this solution include $xH(\frac{x}{y})$ and $xH(\frac{y}{x})$ and $yH(\frac{x}{y})$. All these forms represent homogeneous functions of $x$, $y$ of degree $n=1$. It is easy to see that the solutions $z = \sqrt{xy}$ and $z = x^{\alpha}y^{1-\alpha}$ can be easily obtained by suitable choices of the $H$ function. We also note that the domain of validity $x > 0, y > 0$ required by the $\ln$ function matches that of $z = \sqrt{xy}$ and $z = x^{\alpha}y^{1-\alpha}$.

Finally, we note that for the given domain of validity with $x > 0$, any homogeneous function $f$ of degree $n=1$ can be represented in the form $xH(\frac{y}{x})$. In particular, we have
\begin{gather*}
  f(x, y) = f \left(x\cdot1, x\cdot\frac{y}{x} \right) = xf \left(1, \frac{y}{x} \right) = xH \left(\frac{y}{x} \right)
\end{gather*}
where $H(w) = f(1, w)$.

\section{The Lagrange Integral}

Having looked at Euler's solution, we now look at Lagrange's solution in some detail. The ideas elaborated in this section mirror Goursat's analysis books, \cite[p.426]{Goursat1904} on envelope calculations and \cite[p. 236]{Goursat1917} on details of the complete integral, the general integral, and the singular integral for Clairaut's equation. These are old references written in a different style, so the ideas require a little unpacking. Once we turn those ideas into calculations, we will see that Euler's solution can be obtained from the Lagrange complete integral.

\begin{figure}
  \centering
  \includegraphics[width=0.5 \textwidth]{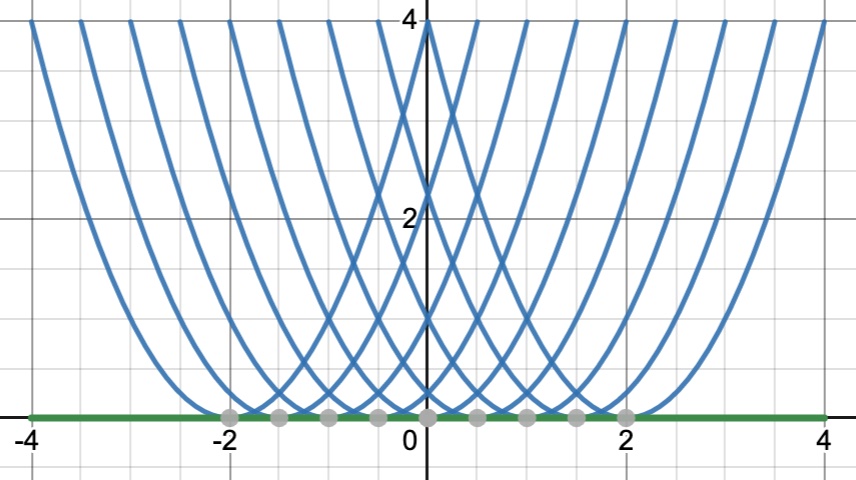}
  \caption{Integrals for $y'^2 = 4y$. The family of parabolas $y = (x + c) ^2$ is called the complete integral, and each member of this family satisfies the PDE. A curve tangent to this family of parabolas is the $x$ axis $y = 0$, and it is called the singular integral, and it satisfies the PDE as well.}
  \label{fig:singular_integral}
\end{figure}

The nature of the singular integral is well explained in \cite{BritannicaIntegral2023} which shows a family of solutions to the equation $y'^2 = 4y$. The family of solutions $y = (x + c)^2$ contains one arbitrary parameter $c$, and can be regarded as a complete integral to the PDE. The tangent to this family of curves is the $x$ axis $y=0$ (\figref{fig:singular_integral}) and it is called the singular integral. It is easy to verify that both the complete integral $y=(x+c)^2$ and the singular integral $y=0$ satisfy the PDE. In this particular example, the complete integral is a higher order curve (a parabola) compared to the singular integral (a straight line). In the case of the restricted Clairaut equation we will find that the complete integral is a simple family of planes while the singular integral takes on rather complex forms based on the parameters chosen.

In the following two subsections, we first derive an envelope calculation procedure, and then use the envelope procedure to calculate singular integrals for the restricted Clairaut equation.

\subsection{A Tangent on Envelopes}

We take a slight digression at this point to describe an envelope calculation procedure. This material typically belongs in a class on Curves and Surfaces, or perhaps in a class on Differential Geometry though we find that the topic is often skipped in the curriculum. We briefly review the relevant procedure and its justification for completeness. More details can be found in older texts like \cite[p.162]{Struik1961} or \cite[p.426]{Goursat1904}.

We would like note here that this article emphasizes clarity over rigor. For instance, the treatment below follows Goursat in describing an envelope procedure. The procedure makes use of certain differentials ($\delta x, \delta y$) which could be expressed more precisely in the language of differential geometry. But Goursat's $\delta$ notation provides good intuition, and we have chosen to follow it without further justification.

\cite[p.426]{Goursat1904} suggests the following. Starting with $f(x, y, a) = 0$, one considers two curves. One curve $C$ that represents $f(x, y, a) = 0$ for a fixed parameter $a$, and the other an envelope $E$ for the entire family of curves $f(x, y, a) = 0$ for different values of $a$. The envelope can be expressed in parameteric form where $x$ and $y$ are functions of $a$, say $x = \phi(a)$ and $y = \psi(a)$. Consider differentials $\delta x$ and $\delta y$ proportional to the directional cosines of the tangent to the curve $C$. In other words, the vector ($\delta x$, $\delta y$) points along the tangent to the curve $C$. But the tangent to the curve $C$ is also tangent to the envelope which has derivatives ($\frac{dx}{da}$, $\frac{dy}{da}$). This gives us the following:
\begin{align*}
    \frac{\frac{dx}{da}}{\delta x} = \frac{\frac{dy}{da}}{\delta y}
\end{align*}

\begin{figure}
  \centering
  \includegraphics[width=0.35 \textwidth]{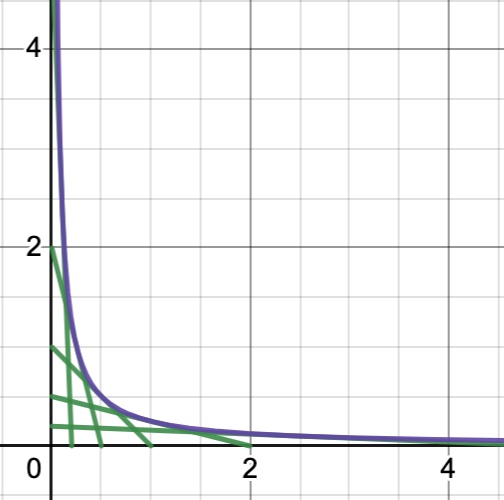}
  \caption{$ax + by = 1$ (with $ab = 1$), $4xy = 1$}
  \label{fig:envelope_example}
\end{figure}

\noindent When we look at the curve $C$, $a$ is fixed and $f$ is constant, so we can infer the following:
\begin{align*}
  df = \frac{\partial f}{\partial x} \delta x + \frac{\partial f}{\partial y} \delta y = 0.
\end{align*}

\noindent Combining the above two we have:
\begin{align} \label{eq:partial_1}
  \frac{\partial f}{\partial x} \frac{dx}{da}  + \frac{\partial f}{\partial y} \frac{dy}{da} = 0.
\end{align}

\noindent This appears to be a neat trick. Using the device of $\delta x$ and $\delta y$, we have obtained an identity that combines properties of the curve ($\frac{\partial f}{\partial x}$, $\frac{\partial f}{\partial y}$) and properties of the envelope ($\frac{dx}{da}$, $\frac{dy}{da}$).

\noindent Now, considering the envelope $E$ directly, $x$ and $y$ are functions of $a$, and so one can write the following:
\begin{align} \label{eq:partial_2}
  \frac{\partial f}{\partial x} \frac{dx}{da} + \frac{\partial f}{\partial y} \frac{dy}{da} + \frac{\partial f}{\partial a} = 0.
\end{align}

\noindent From the above two equations \eqref{eq:partial_1} and \eqref{eq:partial_2} one can infer the envelope equation:
\begin{align} \label{eq:partial_3}
  \frac{\partial f}{\partial a} = 0.
\end{align}

\noindent Further geometric intuition for this envelope equation is provided in \cite[p. 429, Fig. 37]{Goursat1904} which shows nearby curves from the family, and its relation to \eqref{eq:partial_3}. As an example of this procedure, consider the following family of curves:
\begin{align*}
  ax + by = 1.
\end{align*}

\noindent Suppose further that we have a constraint $ab=1$. In this case, we can write:
\begin{align*}
  ax + \frac{1}{a}y = 1, \\
  a^2x + y = a.
\end{align*}

\noindent Now we set:
\begin{align*}
  f(x, y, a) = a^2x + y - a = 0.
\end{align*}

\noindent Letting $\frac{\partial f}{\partial a} = 0$,
\begin{align*}
    2ax - 1 = 0, \\
    a = \frac{1}{2x}.
\end{align*}

\noindent Plugging this back into the family of curves we get the envelope $4xy = 1$ (\figref{fig:envelope_example}):
\begin{align*}
  \frac{1}{2x} \cdot x + 2x \cdot y &= 1, \\
  \text{i.e., } 4xy &= 1.
\end{align*}

\begin{figure}
  \centering
  \includegraphics[width=0.75 \textwidth]{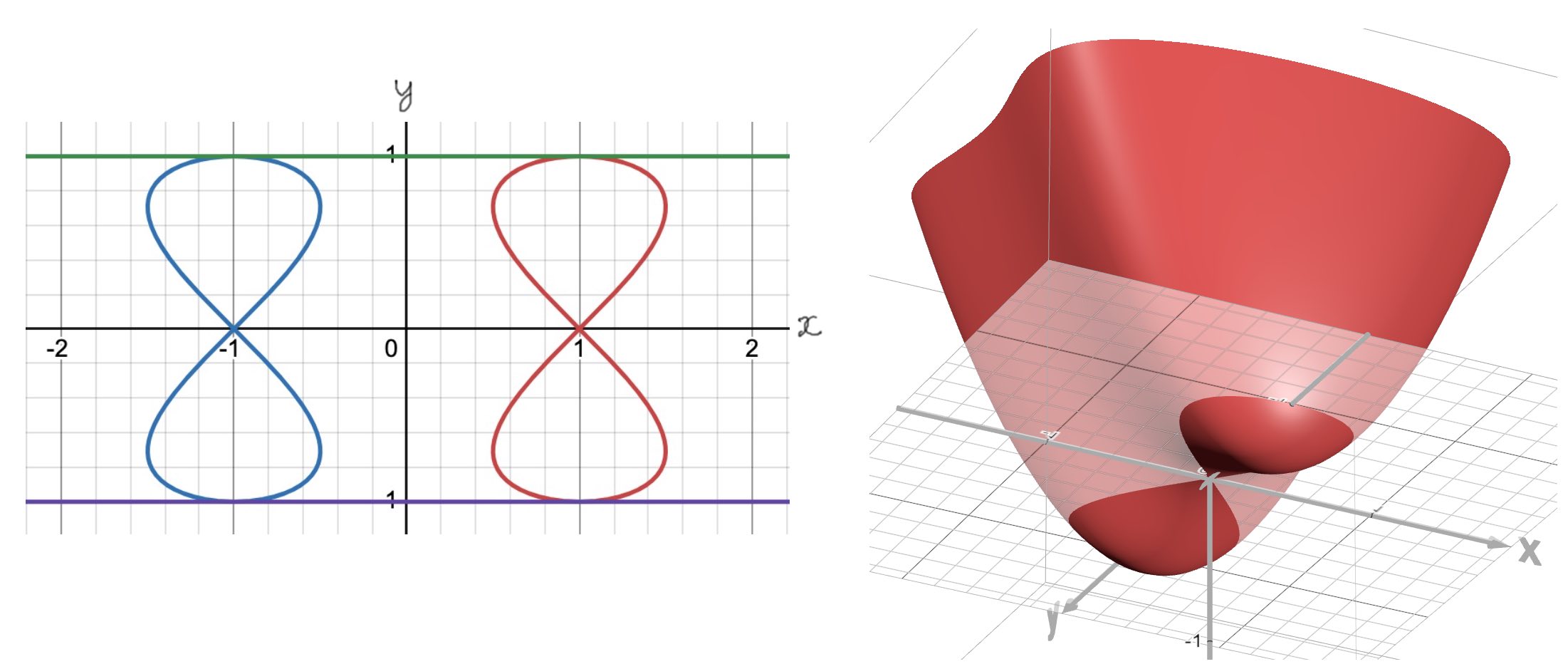}
  \caption{Starting with $y^4 - y^2 - (x - a)^2 = 0$, the line $y=0$ is a locus of singularities, while $y=\pm 1$ represents an envelope}
  \label{fig:goursat_singularity}
\end{figure}

As pointed out in \cite[pp. 427-428]{Goursat1904}, there is some subtlety in associating $\frac{\partial f}{\partial a} = 0$ with an envelope. Applying the condition $\frac{\partial f}{\partial a} = 0$ to $y^4 - y^2 - (x - a)^2 = 0$ we have $y^4 - y^2 = y^2(y^2 - 1) = 0$. This leads to three lines $y = 1$, $y = -1$ and $y = 0$. The first two happen to be envelopes (\figref{fig:goursat_singularity}) while the third is a locus of singularities. The 3-D surface $z = y^4 - y^2$ shows the figure 8 with the saddle point at the origin corresponding to the singularity. We hope this provides some intuition in distinguishing the envelope situation from the singularity situation. In our examples related to the simplified Clairaut equation such loops are not involved, and this type of singularity does not arise.

\subsection{The Singular Integral}

\noindent We now continue our journey to a general solution given the complete integral \eqref{eq:complete_integral} and the envelope equation \eqref{eq:partial_3}. For starters, consider the complete integral $z = ax + by$ where $b$ is a specific function of $a$, say $b = \frac{1}{a}$. This leads to the following:
\begin{align*}
  z = ax + \frac{1}{a}y, \\
  z - ax - \frac{1}{a}y = 0.
\end{align*}

\noindent To get a singular integral, one needs to find the envelope of this single parameter family of solutions:
\begin{align*}
  f(x, y, z, a) = z - ax - \frac{1}{a}y = 0.
\end{align*}

\noindent The envelope is obtained by setting $\frac{\partial f} {\partial a} = 0$. That is:
\begin{align*}
  -x + \frac{1}{a^2}y = 0, \\
  a = \sqrt{\frac{y}{x}}.
\end{align*}

\noindent Using $a = \sqrt{\frac{y}{x}}$, we can now connect the envelope or singular integral back to the earlier solution:
\begin{align*}
  z &= x \sqrt{\frac{y}{x}} + y \sqrt{\frac{x}{y}} = 2 \sqrt{xy}.
\end{align*}

\noindent As a simple exercise, one could try to eliminate the constant 2, that is to obtain $z = \sqrt{xy}$. We could also try to find the solution $z = x^{\alpha}y^{1-\alpha}$ in the above manner.

Now to generalize this discussion, let us consider the complete integral $z = ax + by$ where $b = \phi(a)$ is an arbitrary function of $a$. \cite[p.237]{Goursat1917} provides some justification for this relation between $a$ and $b$ based on which we have the following {\em general} integral:
\begin{align} \label{eq:general_integral}
  z = ax + \phi(a)y.
\end{align}

\noindent Writing this in canonical form,
\begin{align*}
f(x, y, z, a) = z - ax - \phi(a)y = 0.
\end{align*}

\noindent Applying the envelope equation \eqref{eq:partial_3}, namely $\frac{\partial f} {\partial a} = 0$, we obtain the following envelope condition for the above family of planes:
\begin{align} \label{eq:envelope_family}
  x + \phi'(a)y = 0.
\end{align}

\noindent The simplest way to solve \eqref{eq:general_integral} and \eqref{eq:envelope_family} is to write out the solution set in parameteric form as follows:
\begin{align} \label{eq:parameteric_solution}
  (x, y, z) = (-\phi'(a)y, y, -a\phi'(a)y + \phi(a)y).
\end{align}

\noindent Notice that $y$ appears as a factor in all the coordinates, and that all coordinates are zero when $y=0$. So for non-zero $y$, we could divide the parameterized coordinates by $y$, setting $X=\frac{x}{y}$ and $Z=\frac{z}{y}$, to get the following:
\begin{align} \label{eq:projective_solution}
 (X, Z) = (-\phi'(a), -a\phi'(a) + \phi(a)).
\end{align}

\noindent \eqref{eq:projective_solution} is a representation of the solution set in projective space, $\mathbb{R}\mathbb{P}^2$. Though we don't use or define projective spaces in this paper, we mention this representation briefly since it may be a clean and intuitive description of the solution set. Projective representations are useful where homogeneous polynomials are involved with many to many relations $f(x, y, z) = 0$ and we will shortly see solutions of this nature.

\section{Euler-Lagrange Connection} \label{sec:euler_lagrange_connection}

In the first subsection below (Section \ref{sec:euler_from_lagrange}) we derive Euler's homogeneous solution from  Lagrange's complete integral with $b = \phi(a)$ as per the general integral procedure. We also observe that invertible functions $\phi'(a)$ always lead to homogeneous functions as described in Euler's theorem. In Section \ref{sec:convoluted_example} we set $\phi'(a)$ to be a non-invertible function and proceed to construct a different type of solution to the restricted Clairaut equation which is not in Euler's functional form $z=f(x, y)$. These two subsections establish, with some important disclaimers, the first claim, that the Lagrange complete integral is more general than Euler's integral since it can generate Euler's integral, and solutions that go beyond it as well.

The following subsections (Sections \ref{sec:simple_example}, \ref{sec:generalized_general_integral}) contain two simple solutions, in contrast to the convoluted example in Section \ref{sec:convoluted_example}, where simplicity is achieved by relaxing the condition $b=\phi(a)$ in a general integral to other forms like $\phi(a, b) = 0$. This is related to our second claim about extending the general integral, and we further address Evans' comments about the general integral in Section \ref{sec:evans_comments}.

\subsection{Euler's Homogenous Function} \label{sec:euler_from_lagrange}

Let us recall that Euler's homogenous functions satisfy the simplified Clairaut equation. To see this, suppose $z = f(x, y)$ with homogenous $f$ of degree $n=1$, that is, $f(sx, sy) = sf(x, y)$. We differentiate both sides of this equation with respect to $s$ as follows:
\begin{align*}
  \frac{d}{ds} f(sx, sy) &= f_x(sx, sy) \frac{d}{ds}sx + f_y(sx, sy) \frac{d}{ds} sy\\
   &= f_x(sx, sy) x + f_y(sx, sy) y \\
   \frac{d}{ds} sf(x, y) &= f(x, y)
\end{align*}

Setting $s=1$, and equating the right hand sides of the above two equations, we get $f(x, y) = xf_x + yf_y$ as required for the simplified Clairaut equation. In other words, Euler's homogenous functions of degree $1$ satisfy the simplified Clairaut equation.

Given a specific Euler solution we will now deconstruct the corresponding Lagrange integral. This derivation comes with the caveat that certain singularities may limit the domain of the solution as detailed in Section \ref{sec:chojnacki}. Given an Euler integral $z = f(x, y)$ obeying the simplified Clairaut equation, we will prove that the plane $z = \frac{\partial z}{\partial x}|_p x + \frac{\partial z}{\partial y}|_p y$ is tangent to the Euler integral at any given point $p = (x_0, y_0, z_0)$. Firstly, we have $z_0 = \frac{\partial z}{\partial x}|_p x_0 + \frac{\partial z}{\partial y}|_p y_0$ based on the simplified Clairaut equation. Additionally, we have the tangent plane equation $(z - z_0) = \frac{\partial z}{\partial x}|_p (x - x_0) + \frac{\partial z}{\partial y}|_p (y - y_0)$  for a point $(x, y, z)$ on the tangent plane. Combining these two, we have $z = \frac{\partial z}{\partial x}|_p x + \frac{\partial z}{\partial y}|_p y$ for all points $(x, y, z)$ on the tangent plane.

We now have a family of planes  $z = ax + by$ with $a = \frac{\partial z}{\partial x}$ and $b = \frac{\partial z}{\partial y}$ that are tangent to Euler's solution at any given point $(x, y, z)$ as $(x, y)$ range over the full domain $\mathbb{R}^2$. This proves that Euler's solution can be derived from Lagrange's integral in a certain formal sense. We note though an important warning from \cite{Chojnacki1995} that this family of planes may have cusps where the construction breaks down, as in Chojnacki's example $z = \frac{x^3}{y^2}$. His example is of the form $z = xH(\frac{x}{y})$ where $H(\frac{x}{y}) = \frac{x^2}{y^2}$ is a homogenous function as in Euler's solution, and satisfies the simplified Clairaut equation. But as we will see in Section \ref{sec:chojnacki}, Chojnacki describes a problematic cusp at $x=0$ when both partial derivatives vanish and the inverse function theorem no longer applies.

\begin{figure}
  \centering
  \includegraphics[width=0.5 \textwidth]{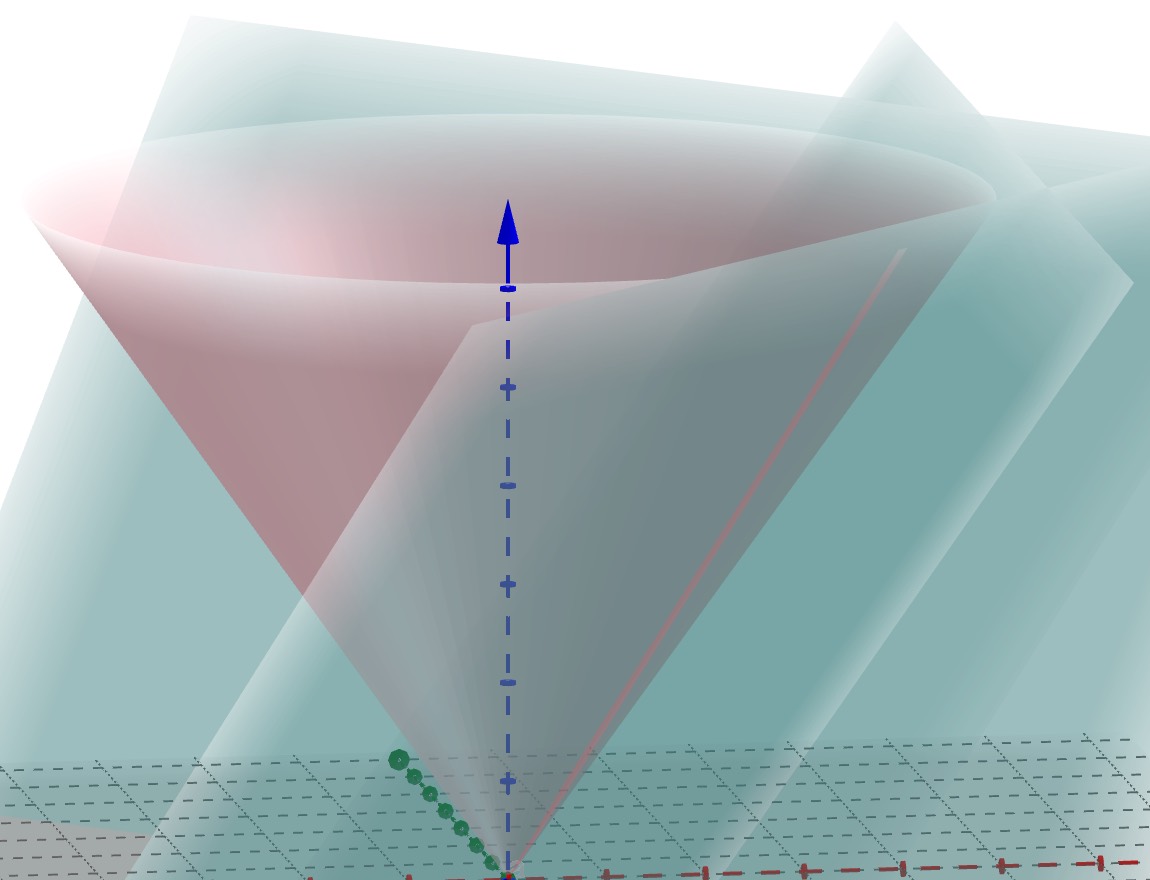}
  \caption{Cone as the envelope of planes}
  \label{fig:cone_envelope}
\end{figure}

We now look at the other direction. Given Lagrange's complete integral we look at the types of solutions it can generate. Let us first recall the general integral \eqref{eq:general_integral}, the envelope of a family of planes for a fixed but arbitrary function $\phi$:
\begin{align} \label{eq:plane_family}
    z = ax + \phi(a) y.
\end{align}

\noindent Note that when both $a$ and $\phi$ are fixed, the above represents the equation of a plane of the form $z = \alpha x + \beta y$ for constants $\alpha, \beta$. Now as we vary $a$ keeping $\phi$ fixed we get a family of planes. It may be helpful to keep \figref{fig:cone_envelope} in mind where three planes from a certain family are shown along with the envelope which is a cone. Also note that points on the envelope satisfy \eqref{eq:envelope_family} and \eqref{eq:plane_family} both of which are homogenous in $(x, y, z)$.

Let the envelope $E$ be tangent to this family of planes at a point $P_1 = (x_1, y_1, z_1)$. Now, $P_1$ obeys \eqref{eq:plane_family} for some fixed plane given by (say) $a=a_1$. It is clear that $sP_1 = (sx_1, sy_1, sz_1)$ lies on the same plane since $sP_1$ also satisfies \eqref{eq:plane_family} with $a=a_1$. Further, when $P_1$ satisfies the envelope condition \eqref{eq:envelope_family}, the point $sP_1$ obeys the envelope condition \eqref{eq:envelope_family} as well, namely $a sx_1 + \phi'(a) sy_1 = 0$. In short, when a point $P=(x, y, z)$ lies on the envelope (satisfying \eqref{eq:envelope_family} and \eqref{eq:plane_family}), the scaled point $sP=(sx, sy, sz)$ lies on the envelope as well, satisfying those two conditions.

If we restrict our solutions to be of the form $z = h(x, y)$, then we also have $sz = h(sx, sy)$ based on the envelope equation for $(sx, sy, sz)$. Multiplying the former by $s$, we have $sz = sh(x, y)$. The two expressions for $sz$ show that $h(sx, sy) = sh(x, y)$. In other $h$ is a homogeneous function, and we have thus derived Euler's solution from the Lagrange integral. Geometrically, this homogeneous solution may be described as a {\em ruled} surface with the radial lines $(sx, sy, sz)$ for fixed $(x, y, z)$ lying on the envelope.

We will now see that invertibility of $\phi'(a)$ leads to Euler's homogeneous functions. We have $\phi'(a) = -\frac{x}{y}$ for non-zero $y$. We have already seen that $y=0$ implies $x=0, z=0$, so assuming non-zero $y$ doesn't lose generality. Now if $\phi'(a)$ is invertible, one could say $a = \psi \left( \frac{x}{y} \right)$ where $\psi$ is the negative inverse of $\phi'$. Plugging this back to the complete integral, we obtain the following singular integral:
\begin{align*}
  z &= \psi \left( \frac{x}{y} \right)x + \phi( \psi \left(\frac{x}{y} \right))y.
\end{align*}

It is clear that the above function is homogeneous of degree $n=1$ in $x$, $y$, and is of the form $z = x H \left( \frac{x}{y} \right)$ for a suitable function $H$. Contrariwise, to go beyond Euler's homogeneous functions we will need to start with $\phi'$ that is not invertible.

\subsection{Beyond Euler's Homogenous Function} \label{sec:convoluted_example}

We have seen that the Lagrange integral yields homogeneous functions when $\phi'(a)$ is invertible. We will now construct a function $\phi'$ that is non-invertible by design, and use it to derive a more general solution to the restricted Clairaut equation. This example is rather convoluted, and we share this primarily for completeness and to illustrate the possibilities and limitations of a general integral based on $b = \phi(a)$. The next section contains simpler algebraic examples.

Let us define $\phi'(a)$ to be the following spline:

\begin{gather}
\phi'(a) =
\left\{
	\begin{array}{ll}
    ((a-1)^2 -1)^2 & \mbox{if } 0 \le a \le 1, \\
    0.5 ((a-1)^2 -1)^2  + 0.5 & \mbox{if } 1 \le a \le 2, \\
    0.5 ((a-3)^2 -1)^2  + 0.5 & \mbox{if } 2 \le a \le 3, \\
    ((a-3)^2 -1)^2 & \mbox{if } 3 \le a \le 4, \\
		0 & \mbox{o.w.}
	\end{array}
\right.
\end{gather}

We have just spliced together four scaled translates of $((a - k)^2 - 1)^2$, each of which is a smooth bump-like function between $k-1$ and $k+1$, for $k=1,3$. This gives us a smooth bimodal function that is non-invertible by design. We obtain $\phi(a)$ from $\phi'(a)$ by a simple numerical integration.

\begin{figure}[H]
  \centering
  \includegraphics[width=0.75\textwidth]{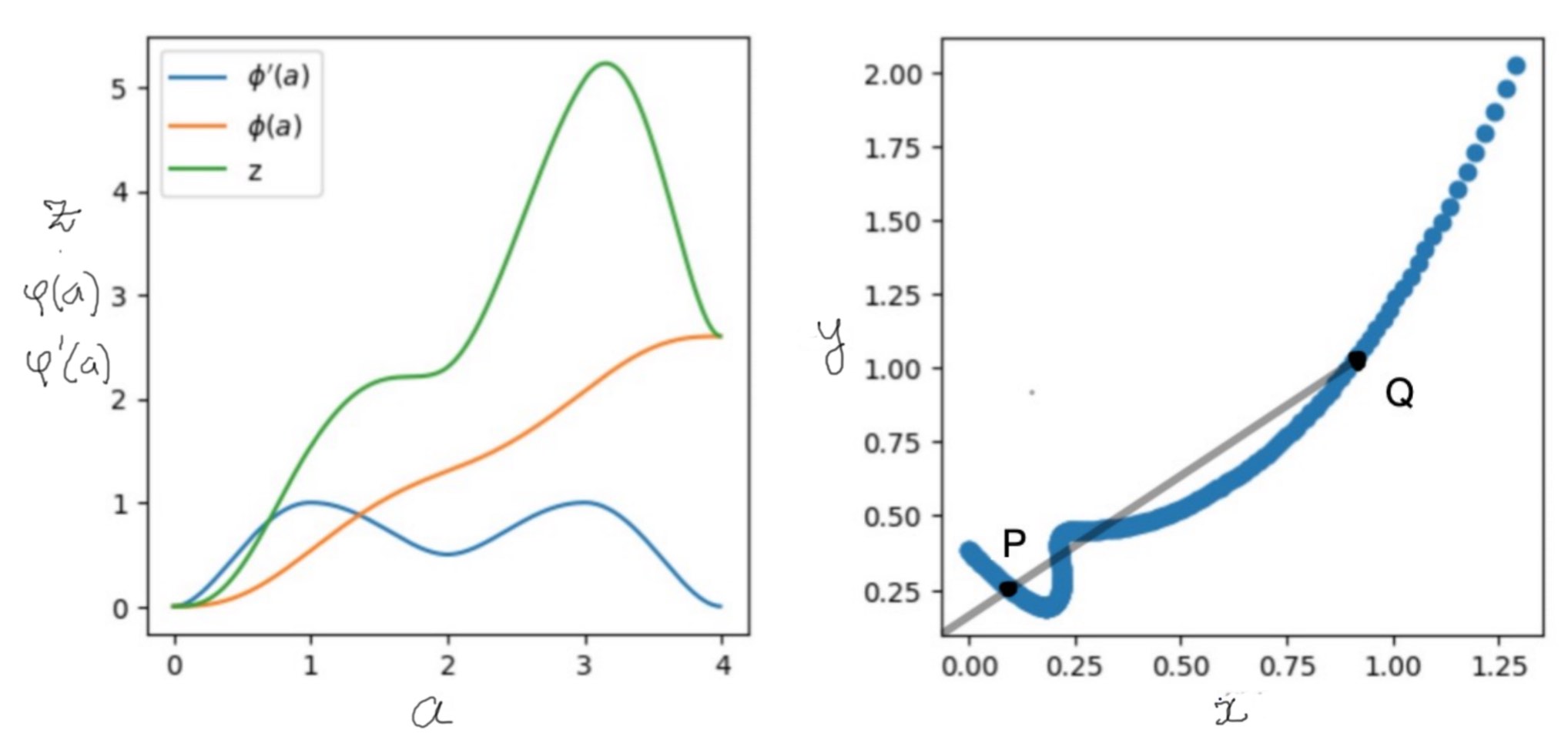}
  \caption{The graphs on the left represent various functions of $a$ starting with an explicitly constructed $\phi'(a)$ that is non-invertible by design. Given $\phi(a)$ as defined on the left, the graph on the right shows a cross section of the surface $z = ax + \phi(a)y$ at $z=1$. $P$ and $Q$ lie on the same radial line, and hence have the same value of $\frac{x}{y}$ in this cross-sectional view. When we project up along the z-axis to get a cone-like surface, this translates to one to many relationship between $(x, y)$ and $z$.}
  \label{fig:clairaut_cross_section}
\end{figure}

Given $\phi(a)$ as above, we seek a cross-sectional view of the surface generated by $z = ax + \phi(a)y$ subject to the envelope condition $\phi'(a) = -\frac{x}{y}$. Holding $y$ fixed at $y=1$, and with $a$ ranging from $0$ to $4$, and $z = ax + \phi(a)y$, we numerically obtain a set of points $(x, y, z) = (-\phi'(a), 1, -a\phi'(a) + \phi(a))$. We then fix $z=1$ and scale the other coordinates suitably, and as per the earlier section, this is valid based on the homogeneity of these solutions. Scaling gives us the following set of points instead $(x, y, z) = (\frac{-\phi'(a)}{-a\phi'(a) + \phi(a)}, \frac{1}{-a\phi'(a) + \phi(a)}, 1)$. With this set of points we plot $x$ against $y$ to get a cross-sectional view at $z=1$ (\figref{fig:clairaut_cross_section}, right).

We now look at two points $P$, $Q$ along a radial line of this cross section as in \figref{fig:clairaut_cross_section}. With $\phi'(a) = -\frac{x}{y}$, the non-invertibility of $\phi'(a)$ guarantees the existence of two such points since radial lines preserve the ratio  $\frac{x}{y}$. We can find two values of $a$ mapping to the same $\phi'(a)$, and thus two values of $(x, y)$ on the same radial line. This is further verified graphically in (\figref{fig:clairaut_cross_section}, right) which was computed numerically. Say $P = (x, y, 1)$ and $Q = (\alpha x, \alpha y, 1)$ with the $z$ coordinate fixed at $1$ for both points. Notice now, that if $P$ is on the surface, so is $sP$, $\forall s \in \mathbb{R}$. Thus $\alpha P = (\alpha x, \alpha y, \alpha)$ is on the surface along with $(\alpha x, \alpha x, 1)$. In particular, we now have two values of $z$, namely $z=1$ and $z=\alpha$ associated with $(\alpha x, \alpha y)$. In other words, $z$ is not a function of $(x, y)$. The solution is a cone-like 3-D surface but it cannot be expressed as a homogeneous function $z=f(x, y)$.

\subsection{A Simple Example} \label{sec:simple_example}

The previous example started with $b=\phi(a)$ and constructed a rather complicated singular integral. Here we start by relaxing the notion of a general integral to include relations $\phi(a, b) = 0$ rather than just $b = \phi(a)$. This allows us to expand our solution set from homogeneous functions to our suggestively named homogeneous surfaces. Let us consider the relation $(a-1)^2 + (b-1)^2 = 1$. To obtain the envelope associated with this relation, we can break this up into two parts each of which is functional. Specifically:
\begin{gather*}
  b = 1 \pm \sqrt{1 - (a-1)^2}
\end{gather*}

\begin{figure}
  \centering
  \includegraphics[width=0.3\textwidth]{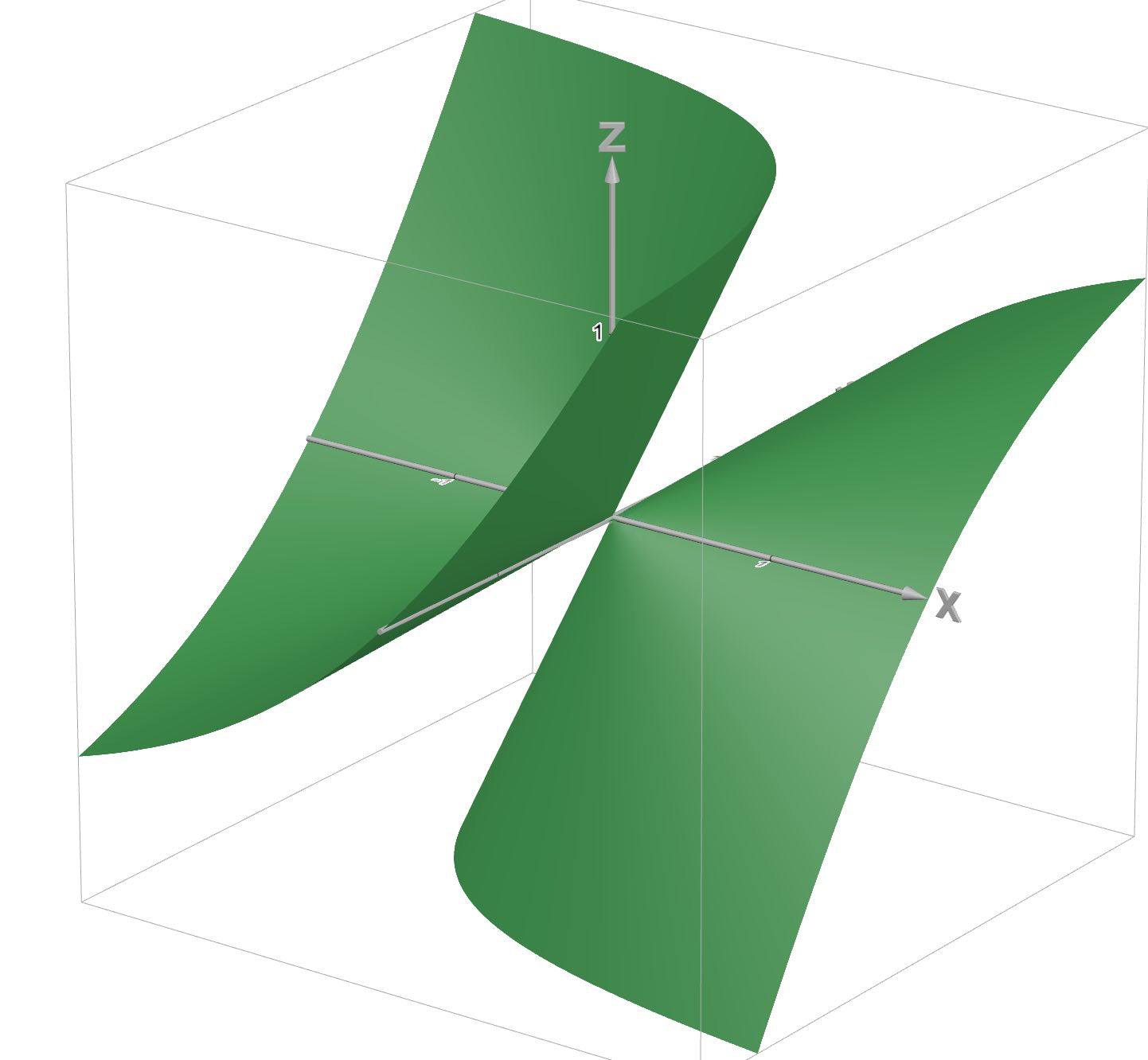}
  \caption{$z^2 = 2xz + 2yz - 2xy$.}
  \label{fig:general_integral_example}
\end{figure}

The two roots represent two functions of $a$. We illustrate our calculation using the positive root since the other one follows the same pattern:
\begin{gather*}
  ax + (1 + \sqrt{1 - (a-1)^2})y = z \\
  \implies f(x, y, a) = ax + (1 + \sqrt{1 - (a-1)^2})y - z = 0
\end{gather*}

Using the envelope equation $\frac{\partial f}{\partial a} = 0$,
\begin{gather} \label{eq:ab_of_xy}
  x + \left( \frac{1}{2 \sqrt{1 - (a-1)^2}} \cdot -2(a-1) \right) y = 0 \\
  a = 1 \pm \frac{x}{\sqrt{x^2 + y^2}}
\end{gather}

We can now plug this back into the complete integral, again with the positive root for simplicity:
\begin{gather*}
  \left( 1 + \frac{x}{\sqrt{x^2 + y^2}} \right)x + \left( 1 + \frac{y}{\sqrt{x^2 + y^2}} \right)y = z \\
  x^2 + y^2 = \sqrt{x^2 + y^2} (z - x - y) \\
  z^2 = 2xz + 2yz - 2xy
\end{gather*}

This surface is shown in \figref{fig:general_integral_example} and accounts for all the positive and negative roots above. We can also verify that the surface $z^2 = 2xz + 2yz - 2xy$ obeys the restricted Clairaut equation as expected. Thus, our expanded notion of the general integral yields singular integrals which go beyond Euler's functional form $z=f(x, y)$.

\subsection{A Generalized General Integral} \label{sec:generalized_general_integral}

\begin{figure}[H]
  \centering
  \includegraphics[width=0.3\textwidth]{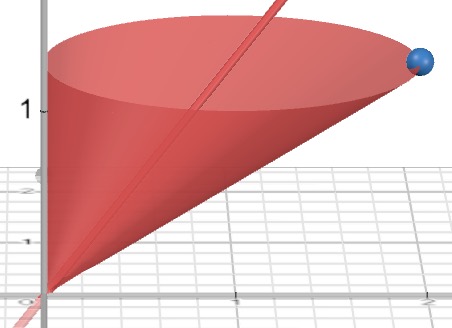}
  \caption{A tilted cone $x^2 + y^2 + z^2 - 2xz -2yz = 0$ with axis $x=y=z$, and cross-sectional radius $r = z$. }
  \label{fig:tilted_cone}
\end{figure}

We now present another solution where we stretch the notion of a general integral beyond the relation $\phi(a, b) = 0$. In this example we start with a singular integral, that is, a solution to the restricted Clairaut equation, and work backwards to find the Lagrange complete integral. A one-line geometric description of this singular integral is that it is a tilted cone, as shown in \figref{fig:tilted_cone}. The axis of the cone is $x = y = z$ and the horizontal $x-y$ cross-sections have radius $r=z$. Algebraically, it can represented in either of these forms:
\begin{gather}
  (x - z)^2 + (y - z)^2 = z^2,  \label{eq:cone_radius_construction}\\
  x^2 + y^2 + z^2 - 2xz - 2yz = 0.
\end{gather}

It is easy to verify algebraically or geometrically that this is not of the functional form $z = f(x, y)$. In particular, the surface has two positive values of $z$ for any given $x, y$ in the first or third quadrants. Let us now verify that this cone satisfies the restricted Clairaut equation. The implicit function theorem tells us that $z$ can be regarded as a function of $x, y$ locally. To find this function we represent the surface as a quadratic equation in $z$:
\begin{gather*}
    z^2 - z(2x + 2y) + (x^2 + y^2) = 0
\end{gather*}

This gives us the below solution which clearly obeys the restricted Clairaut equation since each of the terms $x$, $y$ and $\sqrt{2xy}$ are known to do so:
\begin{gather*}
  z = x + y \pm \sqrt{2xy}
\end{gather*}

Here is another proof using directional derivatives that may provide further intuition. For any such cone, consider a point $P$ along with a horizontal, circular cross-section of the tilted cone. The partial derivatives of $z$ can be understood in terms of radial and tangential coordinates within the cross-sectional plane. Since by construction $z = r$, the cross-sectional radius \eqref{eq:cone_radius_construction}, we can assert the following:
\begin{gather*}
  \frac{\partial z}{\partial r} = 1, \\
  \frac{\partial z}{\partial t} = 0
\end{gather*}

That is, $z$ increases as we move radially (or normally), but remains fixed as we move tangentially along the cross-sectional circle. This radial derivative can be used in the following manner assuming the radius makes an angle $\theta$ with the $x$ axis:
\begin{gather*}
  \frac{\partial z}{\partial x} = \frac{\partial z}{\partial r} \cos{\theta} = \cos{\theta}, \\
  \frac{\partial z}{\partial y} = \frac{\partial z}{\partial r} \sin{\theta} = \sin{\theta}.
\end{gather*}

Plugging these back into the left side of the restricted Clairaut equation:
\begin{align*}
  x \frac{\partial z}{\partial x} + y \frac{\partial z}{\partial y} &= x \cos{\theta} + y \sin{\theta} \\
  &= x \cdot \frac{x}{\sqrt{x^2 + y^2}} + y \cdot \frac{y}{\sqrt{x^2 + y^2}} \\
  &= \sqrt{x^2 + y^2} \\
  &= z
\end{align*}

We have now seen three examples of surfaces not covered by Euler's homogeneous function theorem. Recalling these examples, the first one was based on the notion of a general integral with $b = \phi(a)$ and we proceeded to find a non-invertible $\phi'(a)$ in order to build our example. The next example started with the more general form $\phi(a, b) = 0$. In the current and last example, the family of planes $ax + by = z$ is contrained by some relation between $a$ and $b$ but we have not put it down in explicit form as $\phi(a, b) = 0$. It is easy to see though that the implicit relation between $a$ and $b$ leads to a smooth family of planes. In lieu of a relation $\phi(a, b) = 0$ we can obtain a parameterized version of the tangent plane as follows. Consider a point $P$ on the surface of the cone and two tangent lines at that point. The first tangent $L_1$ is the same as the line $OP$ connecting $P$ with the origin $O$, and it can be represented as follows:
\begin{gather*}
  L_1 =
  \begin{pmatrix}
  t + t \cos{\theta} \\
  t + t \sin{\theta} \\
  t
  \end{pmatrix}
\end{gather*}

The second line $L_2$ is a horizontal tangent at the point $P$ to the circle defined by the horizontal cross-section:
\begin{gather*}
  L_2 =
  \begin{pmatrix}
  t + t \cos{\theta} - \alpha \sin{\theta} \\
  t + t \sin{\theta} + \alpha \cos{\theta} \\
  t
  \end{pmatrix}
\end{gather*}

Since we are looking at tangents to a cone through the origin, the plane defined by $L_1$ and $L_2$ can be equivalently defined in terms of corresponding tangents $L_1'$ and $L_2'$ at a height of $z=1$. In particular, we consider tangents to a point $P'$ where $P'$ lies on $OP$, and has $z$-coordinate $z=1$. The parameterizations of $L_1'$ and $L_2'$ are are a little simpler, and are as follows:

\begin{gather*}
  L_1' =
  \begin{pmatrix}
  1 + \cos{\theta} \\
  1 + \sin{\theta} \\
  1
  \end{pmatrix} \text{ and }
  L_2' =
  \begin{pmatrix}
  1 + \cos{\theta} - \alpha \sin{\theta} \\
  1 + \sin{\theta} + \alpha \cos{\theta} \\
  1
  \end{pmatrix}
\end{gather*}

To parameterize the plane defined by $L_1'$ and $L_2'$, we need the perpendicular $N$ to these two lines which is obtained by computing their cross product:
\begin{gather*}
  N = L_1' \times L_2' =
  \begin{pmatrix}
    (1 + \sin{\theta}) - (1 + \sin{\theta} + \alpha{\cos{\theta}}) \\
    -(1 + \cos{\theta}) + (1 + \cos{\theta} - \alpha{\sin{\theta}}) \\
    (1 + \cos{\theta})(1 + \sin{\theta} + \alpha \cos{\theta}) - (1 + \sin{\theta})(1 + \cos{\theta} - \alpha \sin{\theta})
  \end{pmatrix}
\end{gather*}

After simplification, and dividing throughout by common factor $\alpha$ we have:
\begin{gather*}
  N =
  \begin{pmatrix}
    -\cos{\theta} \\
    -\sin{\theta} \\
    1 + \cos{\theta} + \sin{\theta}
  \end{pmatrix}
\end{gather*}

We can now reduce this to two dimensions by normalizing the $z$ coordinate to $1$ to obtain $N = (a, b, 1)$ where:

\begin{align*}
  a = \frac{-\cos{\theta}}{1 + \cos{\theta} + \sin{\theta}} \\
  b = \frac{-\sin{\theta}}{1 + \cos{\theta} + \sin{\theta}}
\end{align*}

One may notice that the above parameterization blows up at $\theta=\pi$ and $\theta=\frac{-\pi}{2}$. Nevertheless, it is visually clear that the family of planes is smooth, and importantly, the envelope is already known to be the tilted cone. An elementary argument is that $\theta=\pi$ corresponds to the plane $x=0$. In particular, $ax + by = z$ can be rewritten at $x + \frac{b}{a}y = \frac{z}{a}$ which devolves to $x=0$ as $a \to \infty$. Similarly, $\theta = \frac{-\pi}{2}$ corresponds to the plane $y=0$ as $b \to \infty$. But this easy way of explaining away singularities may not always work, or it may require further justification, and pull us towards ideas like local charts. Another approach to establish smoothness of the relationship is to eliminate $\theta$ from the above equations. We refrain from doing so for the sake of argument, as it may be difficult to do in the general case. We now have a smooth parameterized family of planes where the relation between $a$ and $b$ is not shown in the explicit form $\phi(a, b) = 0$, and yet this family of planes yields a valid singular integral, namely the tilted cone.

We now proceed to connect the singular integrals described in this section back to Goursat's general integral. In defining the general integral, Goursat \cite[p.239]{Goursat1917} states that we must choose an arbitrary relation between $a$ and $b$, say $b = \phi(a)$ and proceeds to use this last functional form. It is interesting that Goursat restricted his discussion to the simpler functional form $b = \phi(a)$ for his general integral instead of the more general form $\phi(a, b) = 0$. Our explanation for this is that general relations like $\phi(a, b) = 0$ can be understood locally as a graph $b = \psi(a)$, and Goursat was certainly aware of this given his proof of the inverse function theorem. Further, the use of a general relation like $\phi(a, b) = 0$ introduces difficult questions around the notion of smooth surfaces and smooth families of planes. One of our earlier examples was restricted to a specific relation $(a-1)^2 + (b-1)^2 = 1$ where we avoided this problem by overlaying two functions in a local manner. But speaking more generally of parameterized family of planes, these notions of smoothness may be difficult to define without anticipating the modern theory of manifolds. Goursat may have deftly skirted these problems in his presentation. A general relation between $a$ and $b$ can extend beyond an explicit form such as $\phi(a, b) = 0$ to smooth surfaces patched together like splines or manifolds. We could call the resultant family of planes a generalized general integral and believe that it captures Lagrange's idea more completely.

One may ask if the Lagrange singular integral is the most general one possible, especially in the light of modern differential geometry. This seems like a reasonable proposition for the restricted Clairaut equation, since the Lagrange complete integral is just a smooth collection of tangent planes. One only needs to show that any tangent plane necessarily goes through the origin (so it is of the form $ax + by = z$), but we will not enter into this discussion.

\section{Conclusion}

We went over three singular integral constructions based on the same complete integral $z = ax +  by$. The first was based on the standard general integral formulation $b = \phi(a)$ within which invertible $\phi$ yielded Euler's homogenous function and a more general surface was constructed based on a non-invertible $\phi$. The second construction generalized this to the form $\phi(a, b) = 0$, and the third deconstructed a parameteric form $(a, b) = (f_1(\theta), f_2(\theta))$ from a specific solution $x^2 + y^2 + z^2 - 2xz - 2yz = 0$. These gradual generalizations appear to close the gap between the standard general integral and the complete integral, and partially address the remarks by Evans on the generality of the general integral. There is a separate gap between the space of solutions of a PDE, and what can be expressed via the complete integral. This is the question addressed by Chojnacki, and while our paper acknowledges and notes these limitations, it does not add anything further to this topic.

A lot of care is required when computing envelopes for these different families of planes. We have seen one instance of this in Goursat's example (\figref{fig:goursat_singularity}) where the envelope criterion $\frac{\partial f}{\partial a} = 0$ could describe an envelope, or a locus of singularities. We describe two more tricky situations in the upcoming subsections. The first situation, shown by Chojnacki, is a different type of singularity or cusp where a tangent is not possible due to the nature of the cusp. The second situation, as described by Evans, is an algebraic problem when handling the general integral. We will also connect these arguments to our general integral constructions, and show how we sidestep these problems.

\subsection{Chojnacki on Complete Integrals, Envelopes and Cusps} \label{sec:chojnacki}

Chojnacki \cite{Chojnacki1995} points out that there are difficulties making general statements about the complete integral, and that literature contains some flawed statements in this regard. We discuss Chojnacki's paper briefly to illustrate how these pitfalls have been avoided in our results.

\begin{figure}
  \centering
  \includegraphics[width=0.5 \textwidth]{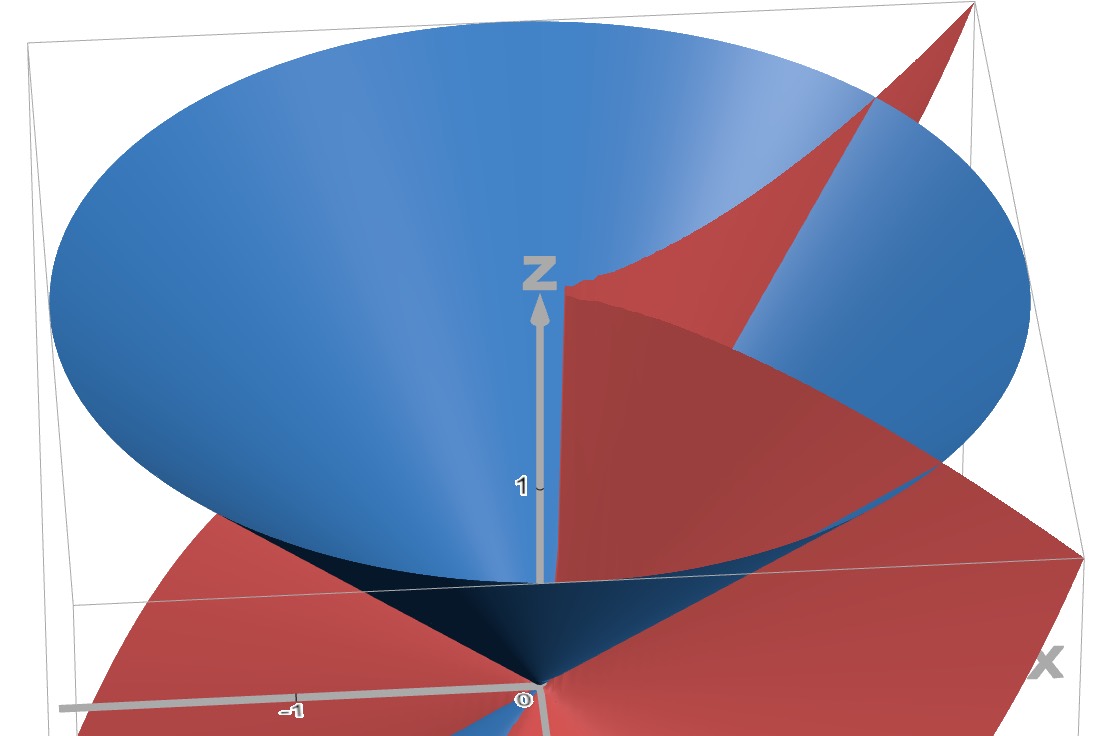}
  \caption{Red surface shows surface $z = \frac{x^3}{y^2}$ with a cusp at $x=0$.}
  \label{fig:chojnacki_cusp}
\end{figure}

Chojnacki's arguments stem from the well known cusp construction $(t^2, t^3)$. The curve is continuous, but not differentiable at $t=0$. This essential construction is generalized to higher dimensions, and the cusp criterion is described in terms of the rank of the differential operator. For our purposes, it suffices to look at the basic two dimensional situation. In this case, Chojnacki starts with a specific solution to the restricted Clairaut equation, namely $z = \frac{x^3}{y^2}$ which is cleverly constructed to have specific partial derivatives and a cusp at $x=0$. Chojnacki's solution is visualized in \figref{fig:chojnacki_cusp}. One can easily imagine that we have a problematic situation describing the tangents at $x=0$, and we elaborate on this point further down.

Reconstructing a family of tangent planes $z = ax + by$ for $z = \frac{x^3}{y^2}$ requires us to specify $(a, b)$ as a function $(x, y)$. Chojnacki recomputes such an inverse function setting $(a, b) = f(x, y) = (f_1(x, y), f_2(x, y)) = (\frac{3x^2}{y^2}, \frac{-2x^3}{y^3})$, the partial derivatives of $z = \frac{x^3}{y^2}$. Substituting these values of $(a, b)$ into $ax + by$ we recover $z = \frac{3x^2}{y^2} \cdot x + \frac{-2x^3}{y^3} \cdot y = \frac{3x^3}{y^2} - \frac{2x^3}{y^2} = \frac{x^3}{y^2}$.

Whereas most of Goursat's analysis regards $(x, y)$ as a function of $(a, b)$, Chojnacki works with the inverse map $f$ which corresponds to $\phi$ in our examples. Such an inverse map $f$ is normally implicit while eliminating $a$, $b$ from Goursat's envelope equations, for instance $a = 1 \pm \frac{x}{\sqrt{x^2 + y^2}}$ \eqref{eq:ab_of_xy} from our earlier example. In Chojnacki's case we have a similar mapping $(a, b) = (\frac{3x^2}{y^2}, \frac{-2x^3}{y^3})$. He then shows that this mapping $(\frac{3x^2}{y^2}, \frac{-2x^3}{y^3})$ (of the form $(3t^2, -2t^3)$) has a cusp at $(0, 0)$ which precludes the possibility of a tangent at that point.

In essence, to use some differential geometry terms, \cite{Chojnacki1995} exhibits a non-differentiable (non-$C^1$) family of tangent planes, and by Frobenius theorem, we can infer that we cannot have a smooth ($C^1$) integral surface. Based on this example \cite{Chojnacki1995} cautions that the notion of tangents and envelopes must be used with care since we may run into such discontinuities. Note that this warning applies to Euler's homogenous functions as well since the solution $z = \frac{x^3}{y^2}$ is of the form $z = xH(\frac{x}{y})$ with $H(\theta) = \theta^2$. Our own subsequent constructions involve continuous families of planes, so these cusp problems don't arise.

One final point about \cite{Chojnacki1995} is that it too constructs a generalized general integral, though in implicit form, and for a different purpose. In particular, Chojnacki constructs a family of planes $z = ax + by$ where $(a, b) = (\frac{3x^2}{y^2}, \frac{-2x^3}{y^3})$. Such a parametric representation $(a, b) = (f_1(x, y), f_2(x, y))$ goes beyond the specific form of a general integral as suggested by \cite{Goursat1917} where $b = \phi(a)$. The relation is not a function, and includes two branches namely $b = \pm 2(\frac{a}{3})^{\frac{3}{2}}$. But Chojnacki provides this construction in passing, and he does not connect it to the notion of a general integral which he doesn't even reference. Further, the existence of such an inverse function $f$ is only asserted point-wise, and locally. Nevertheless it is true that the inverse map $(a, b) = f(x, y)$ is closely connected to the general integral relation $b = \phi(a)$ and to the subsequent generalizations we have described through examples.

\subsection{Evans on the General Integral} \label{sec:evans_comments}

Finally, we note that Evans \cite[p.96, Remarks]{Evans1998} raises some questions about the generality of the general integral which this paper seems to have addressed, at least partially. He provides a theoretical example of a complete integral $F(x, y, z, p, q) = 0$ which could be factorized as
\begin{gather*}
F(x, y, z, p, q) = F_1(x, y, z, p, q) F_2(x, y, z, p, q) = 0.
\end{gather*}
Here $p = \frac{\partial z}{\partial x}$ and $q = \frac{\partial z}{\partial y}$. In such a situation, the general integral for $F_1$ may not cover the solutions of $F_2$, and thus miss certain solutions of the original equation $F(x, y, z, p, q) = 0$.

We find this argument somewhat problematic. If we have a singular integral based on the complete integral for $F_1$, these will indeed miss solutions of $F_2$. But then, there is no reason to think that a complete integral of $F_1$ is also a complete integral of $F$. Here Evans introduces another element, namely the arbitrary function $h$, or to use our terminology $\phi$. But this doesn't add much to the argument since Euler's principle of arbitrary functions doesn't say anything about complete integrals per se. An arbitrary function does not make the complete integral more complete, and we have seen that the principle has some limitations in the first place. It is in fact unclear if the thrust of Evans' argument relates to general integrals, as is Evans' intention, or if it applies to the complete integral as well. To summarize, Evans is asking us to exercise care in making statements about the complete and general integral, a valid point also made by Chojnacki, though more effectively and precisely.

The Goursat general integral procedure could be understood as an algorithm providing local solutions since $b=\phi(a)$ can capture general relationships locally. In the absence of a procedure to patch local solutions together, we need to consider more general relations like $\phi(a, b) = 0$ or other parameteric representations to capture a broader variety of global solutions.

Our constructions of a generalized integral bridge some of the gap between the expressivity of Goursat's general integral, and that of the complete integral. This appears to resolve, at least partially, Evans' primary concern about the generality of the general integral. Our gradual extension of the general integral condition from $b = \phi(a)$ to $\phi(a, b) = 0$ and other parametric forms result in greater expressivity, though we cannot claim to have fully captured the expressivity of the complete integral. A little more work is required by way of proof and construction before we can claim that, and will require some new ideas to go beyond local solutions. Finally, our setup involving Clairaut's equation with $F(x, y, p, q, z) = px + qy - z = 0$ cannot be factorized as $F_1 F_2$ so Evans' specific argument does not apply, at least in a direct way.

\subsection{Final Words}

To summarize the historical development of solutions to the simplified Clairaut equation, we first have Euler's homogeneous function theorem which yields a partial solution. Next comes Lagrange's complete integral $z = ax + by$ and an existence theorem asserting that a general solution to first order PDEs (like the simplified Clairaut equation) can be found by a process of elimination based on the complete integral. Goursat's book describes a specific algorithmic procedure called the general integral (using $b = \phi(a)$) which lets us generate various solutions based on the idea of Lagrange's complete integral. This procedure has made its way to contemporary textbooks like \cite{Evans1998}. Chojnacki notes that it has not been well studied whether Lagrange's complete integral can infact represent all solutions of a given first order PDE, and provides sufficient conditions under which local solutions can be generated using a complete integral. While Chojnacki does not refer to the general integral, Evans notes that the general integral may not be truly general.

With regards to Lagrange's complete integral, the current paper provides a formal derivation of Euler's solution from Lagrange's integral for the simplified Clairaut equation, and also derives further solutions that go beyond Euler's. Despites these arguments and examples, we cannot claim that Lagrange's integral is strictly more general than Euler's. For this, we have Chojnacki's counter example $z = \frac{x^3}{y^2}$ with a cusp at $x=0$. While the Lagrange integral does generate Euler's solution $z = \frac{x^3}{y^2}$ formally, the envelope procedure breaks down in one part of the domain ($x = 0$) because of the cusp. With regards to the general integral, the current paper addresses Evans' arguments and some of these limitations with the general integral procedure. It does this using more general relations like $\phi(a, b) = 0$ and parametric forms like $(a, b) = f(\theta)$ or $(a, b) = f(x, y)$, thus bridging the gap between the general integral and the complete integral.

Our arguments can also be understood in the vocabulary of differential geometry. The complete integral $z = ax + by$ with suitably varying $a, b$ corresponds to a family of tangent planes, in effect a distribution. An envelope or singular integral for this set of tangent planes is analogous to an integral surface or submanifold, as in Frobenius' integrability theorem. Such a viewpoint may provide a good pedagogical introduction to Frobenius' theorem, and seems like an interesting direction for future work.

\begin{acknowledgement}
We would like to thank Prof. Rajaram Nithyananda for suggesting the elegant $X = \ln x$ solution and Dr. Gobinda Sau for detailed feedback on the paper, and for help deconstructing \cite{Chojnacki1995}.
\end{acknowledgement}

\printbibliography
\end{document}